\documentstyle[12pt]{article}

 \newcommand{\sect}[1]{\setcounter{equation}{0}\section{#1}}

 \font\tengoth=eufm10 \font\sevengoth=eufm7 \font\fivegoth=eufm5
 \newfam\gothfam
 \textfont\gothfam=\tengoth \scriptfont\gothfam=\sevengoth
   \scriptscriptfont\gothfam=\fivegoth
   \def\goth{\fam\gothfam}    
 \font\frak=eufm10 scaled\magstep1
 \newcommand{\fra}[1]{\mbox{\frak #1}}

 \def\be{\begin{equation}}
 \def\ee{\end{equation}}
 \def\bea{\begin{eqnarray}}
 \def\eea{\end{eqnarray}}

\def\g{\mbox{${\fra g}$}}

\def\gz{\mbox{${\fra g}_{0}$ }}
\def\gzcar{\mbox{${\fra g}_{0}^{\rm car}$ }}
\def\gkcar{\mbox{${\fra g}_{k}^{\rm car}$ }}

\def\caz{\mbox{${\cal A}_{0}$ }}
\def\cazcar{\mbox{${\cal A}_{0}^{\rm car}$ }}
\def\cao{\mbox{${\cal A}_{1}$ }}
\def\cakp{\mbox{${\cal A}_{k+1}$ }}
\def\cak{\mbox{${\cal A}_{k}$ }}
\def\cakcar{\mbox{${\cal A}_{k}^{\rm car}$ }}
\def\cbk{\mbox{${\cal B}_{k}$ }}
\def\cbz{\mbox{${\cal B}_{0}$ }}
\def\cbo{\mbox{${\cal B}_{1}$ }}
\def\cbkp{\mbox{${\cal B}_{k+1}$ }}

\def\hlamnk{\mbox{${H_{\lambda_{0}^k}}$ }}
\def\hlamnkm{\mbox{${H_{\lambda_{0}^{k-1}}}$ }}
\def\hlamnz{\mbox{${H_{\lambda_{0}^{0}}}$ }}

 \def\L{\mbox{${\bf L}$}}

\def\eiphi{\mbox{$\Phi_{{\cal E}_i}$}}
\def\ekphi{\mbox{$\Phi_{{\cal E}_k}$}}
\def\ekmphi{\mbox{$\Phi_{{\cal E}_{k-1}}$}}
\def\ezphi{\mbox{$\Phi_{{\cal E}_0}$}}
\def\eophi{\mbox{$\Phi_{{\cal E}_1}$}}
\def\ephi{\mbox{$\Phi_{\cal E}$}}
\def\jphi{\mbox{$\Phi_{\cal J}$}}
\def\jkphi{\mbox{$\Phi_{{\cal J}_{k}}$}}
\def\jzphi{\mbox{$\Phi_{{\cal J}_{0}}$}}
\def\jophi{\mbox{$\Phi_{{\cal J}_{1}}$}}

\def\fae{\mbox{${\cal F}_{\cal E}$ }}
\def\faethe{\mbox{${\cal F}_{{\cal E}_{\Theta}}$ }}

\def\fabk{\mbox{${\cal F}_{{\cal B}_{k}}$}}

\def\fabkp{\mbox{${\cal F}_{{\cal B}_{k+1}}$}}
\def\fabkm{\mbox{${\cal F}_{{\cal B}_{k-1}}$}}

\def\chfabk{\mbox{${\cal F}_{{\cal B}_{k \prec 0}}$ }}
\def\chfabo{\mbox{${\cal F}_{{\cal B}_{1 \prec 0}}$ }}
\def\chfabkm{\mbox{${\cal F}_{{\cal B}_{(k-1) \prec 0}}$ }}
\def\chfabp{\mbox{${\cal F}_{{\cal B}_{p \prec 0}}$ }}

\def\faeab{\mbox{${\cal F}_{\cal E(\alpha,\beta)}$}}
\def\faesab{\mbox{${\cal F}_{\cal E'(\alpha,\beta)}$ }}

\def\cherbk{\mbox{${\cal R}_{{\cal B}_{k \prec 0}}$ }}
\def\cherbp{\mbox{${\cal R}_{{\cal B}_{p \prec 0}}$ }}

\def\chercbp{\mbox{$ r_{{\cal B}_{p \prec 0}}$ }}

\def\lamnk{\mbox{$\lambda_{0}^k$ }}
\def\lamnkm{\mbox{$\lambda_{0}^{k-1}$ }}
\def\lamnz{\mbox{$\lambda_{0}^0$ }}
\def\lamno{\mbox{$\lambda_{0}^1$ }}
\def\signk{\mbox{$\sigma_{0}^k$ }}
\def\signkm{\mbox{$\sigma_{0}^{k-1}$ }}
\def\signz{\mbox{$\sigma_{0}^0$ }}
\def\kpi{\mbox{$\pi_k$ }}

\def\kmpip{\mbox{$\pi'_{k-1}$ }}
\def\zpi{\mbox{$\pi_0$ }}
\def\zpip{\mbox{$\pi'_0$ }}
\def\opi{\mbox{$\pi_1$ }}

 \def\deleabh{\mbox{$\Delta_{\cal E(\alpha , \beta)}\,(H)$}}
 \def\deleaba{\mbox{$\Delta_{\cal E(\alpha , \beta)}\,(A)$}}
 \def\deleabb{\mbox{$\Delta_{\cal E(\alpha , \beta)}\,(B)$}}
 \def\deleabe{\mbox{$\Delta_{\cal E(\alpha , \beta)}\,(E)$}}

 \def\delesabh{\mbox{$\Delta_{\cal E'(\alpha , \beta)}\,(H)$}}
 \def\delesaba{\mbox{$\Delta_{\cal E'(\alpha , \beta)}\,(A)$}}
 \def\delesabb{\mbox{$\Delta_{\cal E'(\alpha , \beta)}\,(B)$}}
 \def\delesabe{\mbox{$\Delta_{\cal E'(\alpha , \beta)}\,(E)$}}

\begin{document}
\begin{center}

{\LARGE{\bf{Chains of twists for classical Lie algebras}}}
\footnote{This  work
has been partially supported by DGES of the  Ministerio de  Educaci\'on y
Cultura of  Espa\~na under Projects PB95-0719 and  SAB1995-0610,
the Junta de
Castilla y Le\'on (Espa\~na), and the  Russian  Foundation for Fundamental
Research under grants  97-01-01152 and 98-01-00310.}\\[3mm]

\vskip1cm

{\sc Petr P. Kulish }
\vskip0.25cm

{\it  St. Petersburg Department of the Steklov Mathematical
Institute,}\\
{\it 191011, St.Petersburg, Russia}
\vskip0.5cm

{\sc Vladimir D. Lyakhovsky}
\vskip0.25cm
{\it Theoretical Department, St. Petersburg State University,}\\
{\it 198904, St. Petersburg, Russia.}
\vskip0.5cm

{\sc Mariano A. del Olmo}
\vskip0.25cm
{\it  Departamento de  F\'{\i}sica Te\'orica,
Universidad de Valladolid,  }\\
{\it E-47011, Valladolid,  Spain}

\vskip0.6cm

{emails: kulish@pdmi.ras.ru, lyakhovs@snoopy.phys.spbu.ru,
olmo@fta.uva.es}
 \end{center}

\vskip1cm
\centerline{\today}
\vskip1.5cm

\begin{abstract}
For chains of regular injections ${\cal A}_p \subset {\cal A}_{p-1}
\subset \ldots \subset {\cal A}_{1} \subset {\cal A}_{0}$
of Hopf algebras the sets of maximal
extended Jordanian twists $\left\{ {\cal F}_{{\cal E}_k} \right\}$
are considered.
We prove that under certain conditions there exists for ${\cal A}_{0}$
the twist  \chfabk composed by the factors ${\cal F}_{{\cal E}_k}$.
The general construction of a chain of twists is applied to the
universal envelopings $U(\g)$ of classical Lie algebras $\g$ . We study the
chains for the infinite series $A_n ,B_n$ and $D_n$. The properties
of the deformation produced by a chain $U(\g)_{\scriptsize{\chfabk}}$
are explicitly demonstrated for the case of $\g = so(9)$.
\end{abstract}
\newpage

\sect{Introduction}

The triangular Hopf algebras and twists (they preserve the triangularity
\cite{D2,D3}) play an important role in quantum group theory and
applications \cite{KUL1,VLA,VAL}. Very few types of twists were written
explicitly in a closed form. The well known example is the jordanian twist
(${\cal JT}$) of the Borel algebra $B(2)$ ($\{H,E|[H,E]=E\}$)
with $r=H\wedge E$ \cite{DRIN} where the triangular
$R$--matrix ${\cal R}=(\jphi)_{21}\jphi^{-1}$  is defined by the twisting
element \cite{OGIEV,GER}
\begin{equation}
\label{og-twist}
\jphi=\exp \{H\otimes \sigma \},
\end{equation}
with $\sigma = \ln (1 + E)$.
In \cite{KLM} it was shown that there
exist different extensions (${\cal ET}$'s)
of this twist.
Using the notion of factorizable twist \cite{RSTS}
the element $\fae \in {\cal U} (sl(N))^{\otimes 2}$,
\begin{equation}
\label{twist-sl(N)}
\fae= \Phi_{\cal E} \jphi = \left( \prod^{N-2}_{i=2}\eiphi \right) \jphi
=\exp \{2\xi \sum_{i=2}^{N-1}E_{1i}\otimes
E_{iN}e^{-\widetilde{\sigma} }\}\exp \{H\otimes \widetilde{\sigma} \},
\end{equation}
was proved to satisfy the twist equation,
where $E=E_{1N}$, $H=E_{11}-E_{NN}$ is one of the Cartan
generators $ H \in {\goth h}(sl(N))$,
$\widetilde{\sigma} =\frac 12\ln (1+2\xi E)$ and
$\{ E_{ij} \} _{i,j = 1, \dots ,N} $ is the standard $gl(N)$ basis.

Any subset of exponentials $\left\{ \eiphi | i=2,\dots,N-2 \right\}$
can be used to form an extended twist like (\ref{twist-sl(N)}).
This means that similar extended twistings can be applied to different
algebras with similar structure. In this particular case it is not
difficult to explain this effect: the factors \eiphi \ commute and
the subalgebras where they are defined (the carrier subalgebras
\cite{GER}) intersect by the central element $E$.

Let ${\cal A}$ be a Hopf algebra, ${\cal B}$ and ${\cal C}$ be its
subalgebras such that they are carriers for twists ${\cal F_B}$ and  ${\cal
F_C}$ respectively.
It is important to know under what conditions the sequence ${\cal
F_C F_B}$
provides a new twisting element and what are its properties.
In this paper we study the possibility to compose extended twists
for the universal enveloping algebras of classical Lie algebras.

In Section 2 we present a short list of basic relations for twists.
The general properties of extended twists  are displayed
in Section 3. The sufficient conditions for the existence of a
composition of twists defined for subalgebras are formulated in Section 4.
In Section 5 the same problem is solved for chains of subalgebras.
It is proved that the corresponding chains of twists \chfabk exist
in classical Lie algebras of the series $A,B$ and $D$. Using the regular
injection $A_{n-1} \longrightarrow C_n$ one can implement into $U(C_n)$
a chain typical for $A_{n-1}$. Such improper chains are studied in Section 6.
The properties of twisting performed by a chain \chfabk  are illustrated by
the explicit example of a deformation
$U(so(9)) \stackrel{\scriptsize{\chfabk}}{\longrightarrow}
U(so(9))_{\scriptsize{\chfabk}} $
presented in Section 7. We conclude with some brief remarks about the
possible multiparametrization of chains and the corresponding deformations.


\sect{Basic definitions }

In this section we remind briefly the basic properties of twists.

A Hopf algebra ${\cal A}(m,\Delta ,\epsilon,S)$ with
multiplication $m\colon {\cal A}\otimes {\cal A}\to {\cal A}$,
coproduct $\Delta \colon {\cal A}\to {\cal A}%
\otimes {\cal A}$, counit $\epsilon \colon {\cal A}\to C$,
and antipode $S : {\cal A}\to {\cal A}$
can be transformed \cite{D2} by an invertible (twisting)
element ${\cal F}\in {\cal A}
\otimes {\cal A}$, ${\cal F}=\sum f_i^{(1)}\otimes f_i^{(2)}$,
into a twisted
one ${\cal A}_{\cal F}(m,\Delta _{\cal F},\epsilon ,S_{\cal F})$.
This Hopf algebra ${\cal A}_{\cal F}$ has the
same multiplication and counit  but the twisted coproduct and antipode given by
\begin{equation}
\label{def-t}
\Delta _{\cal F}(a)={\cal F}\Delta (a){\cal F}^{-1},\qquad S_{\cal
F}(a)=vS(a)v^{-1},
\end{equation}
with
$$
v=\sum f_i^{(1)}S(f_i^{(2)}), \qquad a\in {\cal A}.
$$
The twisting element has to satisfy the equations
\begin{eqnarray}
\label{def-n}
(\epsilon \otimes  id)({\cal F}) = (id \otimes  \epsilon)({\cal F})=1,
\\[0.2cm]
\label{gentwist}
{\cal F}_{12}(\Delta \otimes  id)({\cal F}) =
{\cal F}_{23}(id \otimes  \Delta)({\cal F}).
\label{TE}
\end{eqnarray}
The first one is just a normalization condition and
follows from the second relation modulo a non-zero scalar factor.

If ${\cal A}$ is a Hopf subalgebra of ${\cal B}$ the twisting
element ${\cal F}$
satisfying (\ref{def-n}) and (\ref{gentwist}) induces the twist
deformation  ${\cal B}_{\cal F}$ of  ${\cal B}$. In this case one can
put $a \in  {\cal B}$ in all the formulas (\ref{def-t}). This will
completely define the Hopf algebra ${\cal B}_{\cal F}$. Let ${\cal A}$ and
${\cal B}$ be the universal enveloping algebras: ${\cal A} = U({\fra l})
\subset {\cal B}= U(\g)$ with ${\fra l} \subset \g$. If $U({\fra l})$
is the minimal subalgebra on which ${\cal F}$ is completely defined
as ${\cal F} \in U({\fra l})
\otimes U({\fra l})$ then ${\fra l}$ is called the carrier algebra for
${\cal F}$ \cite{GER}.

The composition of appropriate twists can be defined as
${\cal F} = {\cal F}_2 {\cal F}_1$. Here the element ${\cal F}_1$ has to
satisfy the twist equation with the coproduct of the original Hopf algebra,
while ${\cal F}_2$ must be its solution for $\Delta_{{\cal F}_1}$ of the
algebra twisted by ${\cal F}_1$.

If the initial Hopf algebra ${\cal A}$ is quasitriangular with the
universal element ${\cal R}$ then so is the twisted one
${\cal A}_{\cal F}(m,\Delta _{\cal F},\epsilon,S_{\cal F},{\cal R}_{\cal F})$
with
\begin{eqnarray}\label{Rt}
{\cal R}_{\cal F}=  {\cal F}_{21} \,{\cal R} \,{\cal F}^{-1}.
\end{eqnarray}

Most of the explicitly known twisting elements have the factorization
property with respect to comultiplication
$$
(\Delta \otimes id)({\cal F})={\cal F}_{23}{\cal F}_{13}\qquad \mbox{or}
\qquad
(\Delta \otimes id)({\cal F})={\cal F}_{13}{\cal F}_{23}\,,
$$
and
$$
(id \otimes \Delta)({\cal F})={\cal F}_{12}{\cal F}_{13}\qquad \mbox{or}
\qquad
(id \otimes \Delta)({\cal F})={\cal F}_{13}{\cal F}_{12}\,.
$$
To guarantee the validity of the twist equation, these identities are to be
combined with  the additional requirement
${\cal F}_{12}{\cal F}_{23}={\cal F}_{23}{\cal F}_{12}$
or the Yang--Baxter equation on ${\cal F}$ \cite{RSTS,RES}.

An important subclass of factorizable twists consists of elements
satisfying the equations
\begin{eqnarray} \label{f-twist1}
(\Delta \otimes id)({\cal F})={\cal F}_{13}{\cal F}_{23}\,,
\\ [0,2cm] \label{f-twist2}
(id\otimes \Delta _{\cal F})({\cal F})={\cal F}_{12}{\cal F}_{13 }\,.
\end{eqnarray}
Apart from the universal $R$--matrix ${\cal R}$ that satisfies these
equations for $\Delta_{\cal F}=\Delta ^{op}$ ($\Delta ^{op}=\tau\circ \Delta$,
where $\tau(a\otimes b)=b\otimes a$)  there are two
more well developed  cases of such twists: the jordanian twist of the Borel
algebra $B(2)$   where ${\cal F}_j$
has the form (\ref{og-twist}) (see \cite{OGIEV}) and the extended
jordanian  twists (see \cite{KLM} and \cite{V-M,KM} for details).

According to the result by Drinfeld \cite{D3}  skew (constant)
solutions of the classical Yang--Baxter equation (CYBE) can be quantized and
the  deformed algebras thus  obtained can be presented in a form of twisted
universal enveloping  algebras. On the other hand, such solutions of CYBE
can be connected with  the quasi-Frobenius carrier subalgebras of the initial
classical Lie  algebra \cite{STO}. A Lie  algebra $\g(\mu)$,
with the Lie composition $\mu$,
is called Frobenius if there exists a linear functional $g^* \in \g^*$
such that the form $b(g_1,g_2)=g^*(\mu(g_1,g_2))$ is nondegenerate.
This means that $\g$ must have a nondegenerate 2--coboundary $b(g_1,g_2) \in
B^2(\g,{\bf K})$. The algebra is called quasi-Frobenius if it has a
nondegenerate 2--cocycle $b(g_1,g_2) \in Z^2(\g,{\bf K})$ (not
necessarily a coboundary). The classification of quasi-Frobenius
subalgebras in $sl(n)$ was given in \cite{STO}.

The deformations of quantized algebras include the deformations of their
Lie bialgebras $(\g,\g^*)$. Let us fix one of the constituents
$\g^*_1(\mu^*_1)$ (with composition $\mu^*_1$) and deform it  in the  first
order
$$
(\mu^*_1)_t = \mu^*_1 + t\mu^*_2,
$$
its deforming function $\mu^*_2$ is also a Lie product and the deformation
property becomes reciprocal: $\mu^*_1$ can be considered as a first order
deforming function  for the algebra $\g^*_2(\mu^*_2)$. Let $\g(\mu)$
be a Lie algebra that form Lie bialgebras with  both
$\g^*_1$ and $\g^*_2$. This means that we have a one-dimensional
family $\{ (\g,(\g^*_1)_t) \}$ of Lie bialgebras and correspondingly a
one-dimensional family of quantum deformations
$ \{{\cal A}_t(\g,(\g^*_1)_t) \}$
\cite{ETI}. This situation provides the possibility to construct in the
set of Hopf algebras a smooth curve connecting quantizations of the
type ${\cal A}(\g,\g^*_1)$ with those of ${\cal A}(\g,\g^*_2)$. Such  smooth
transitions can involve contractions provided
$\mu^*_2 \in B^2 (\g^*_1,\g^*_1)$. This happens in the case of
${\cal JT,\ ET}$ and some other twists (see \cite{KL} and references therein).

\sect{Extended twists}

Extended jordanian twists are associated with the set
$\{{\bf L}(\alpha,\beta,\gamma,\delta)_{\alpha + \beta = \delta}\}$
of Frobenius algebras \cite{KLM},\cite{V-M}
\begin{equation}
\begin{array}{l}
 [H,E] = \delta E, \quad [H,A] = \alpha A, \quad [H,B] = \beta B, \\[0.2cm]
[A,B] = \gamma E, \quad[E,A] = [E,B] = 0, \quad \quad
\alpha + \beta = \delta .
\end{array}
\end{equation}
For limit values of $\gamma$ and $\delta$ the structure of ${\bf L}$ 
degenerates. For the internal (nonzero) values of $\gamma$ and $\delta$ the
twists associated with the corresponding ${\bf L}$'s are equivalent. It is 
sufficient to study the one-dimensional subvariety 
${\cal L}=\{{\bf L}(\alpha,\beta)_{\alpha + \beta = 1}\}$,
that is to consider the carrier algebras 
\begin{equation}
\label{l-norm}
\begin{array}{l}
 [H,E] = E, \quad [H,A] = \alpha A, \quad [H,B] = \beta B, \\[0.2cm]
[A,B] = E, \quad[E,A] = [E,B] = 0, \quad \quad
\alpha + \beta = 1 .
\end{array}
\end{equation}
The corresponding group 2--cocycles (twists) are
\be
\label{t-ext}
\faeab =  \Phi_{\cal E(\alpha, \beta)} \Phi_j
\ee  
or
\be
\label{t-ext-s}
\faesab = \Phi_{\cal E'(\alpha, \beta)}  \Phi_j
\ee  
with
\begin{equation}
\label{fractions}
\begin{array}{lcl}
 \Phi_j & = &  \exp \{H\otimes \sigma \}, \\[0.2cm]
 \Phi_{\cal E(\alpha, \beta)} & = &
\exp \{ A \otimes B e^{-\beta \sigma} \}, \\[0.2cm]
 \Phi_{\cal E'(\alpha, \beta)} & = &
\exp \{ -B \otimes A e^{-\alpha \sigma} \}.
\end{array}
\end{equation}
Twists (\ref{t-ext}) and (\ref{t-ext-s})
 define the deformed Hopf algebras ${\bf L}_{\cal E(\alpha, \beta)} $
with the co-structure
\begin{equation}
\label{e-costr}
\begin{array}{lcl}
 \deleabh & = & H \otimes e^{-\sigma} + 1 \otimes H
              - A \otimes B e^{-(\beta + 1)\sigma}, \\[0.2cm]
 \deleaba & = & A \otimes e^{-\beta \sigma} + 1 \otimes A , \\[0.2cm]
 \deleabb & = & B \otimes e^{\beta \sigma} + e^{\sigma} \otimes B , \\[0.2cm]
 \deleabe & = & E \otimes e^{\sigma} + 1 \otimes E =
                E \otimes 1 + 1 \otimes E + E \otimes E;
\end{array}
\end{equation}
 and
$ {\bf L}_{\cal E'(\alpha, \beta)} $ defined by
\begin{equation}
\label{e-pr-costr}
\begin{array}{lcl}
 \delesabh & = & H \otimes e^{-\sigma} + 1 \otimes H
              + B \otimes A e^{-(\alpha + 1)\sigma}, \\[0.2cm]
 \delesaba & = & A \otimes e^{\alpha \sigma} + e^{\sigma} \otimes A ,
\\[0.2cm]
 \delesabb & = & B \otimes e^{-\alpha \sigma} + 1 \otimes B , \\[0.2cm]
 \delesabe & = & E \otimes e^{\sigma} + 1 \otimes E =
                  E \otimes 1 + 1 \otimes E + E \otimes E.
\end{array}
\end{equation}
The sets  $\{ {\bf L}_{\cal E(\alpha, \beta)} \}$ and
$\{ {\bf L}_{\cal E'(\alpha, \beta)} \}$
 are equivalent due to the Hopf isomorphism
${\bf L}_{\cal E(\alpha, \beta)} \approx {\bf L}_{\cal E'(\beta, \alpha)}$:
\be
\{ {\bf L}_{\cal E}(\alpha,\beta) \} \approx \{ {\bf L}_{\cal E'}
(\alpha,\beta) \}
\approx \{ {\bf L}_{\cal E}(\alpha \geq \beta) \} \cup \{ {\bf L}_{\cal E'}
(\alpha \geq \beta) \}.
\ee
So, it is sufficient to use only one of the extensions
either $\Phi_{\cal E(\alpha, \beta)}$ or $\Phi_{\cal E'(\alpha, \beta)}$,
or a half of the domain for $(\alpha, \beta)$.

The set ${\cal L}=\{{\bf L}(\alpha,\beta)_{\alpha + \beta = 1}\}$  is just
the family of 4-dimensional Frobenius algebras that one finds in $U(sl(N))$
\cite{STO}.

\sect{Sequences of twists}

Consider again the formula (\ref{twist-sl(N)}) (now on we use a basis
normalized as in (\ref{l-norm}), so  here $H= \frac{1}{2} (E_{11} - E_{NN})$
),
\begin{equation}
\label{twist-sl(N)-fac}
\fae= \left( \prod_{i=2}^{N-1}\eiphi \right) \jphi =
\left( \prod_{i=2}^{N-1}\exp
\{ E_{1i}\otimes E_{iN}e^{-\frac{1}{2}{\sigma} }\} \right) \exp \{H\otimes
{\sigma}
\}.
\end{equation}
In the product of exponentials each factor $\Phi_{{\cal E}_r}$ is itself a
twisting element for the Hopf algebra previously twisted by $\left(
\prod_{i=2}^{r-1}\eiphi \right) \jphi$. This is a very simple example  of a
chain of twists. All the factors \eiphi \ commute
and the corresponding twistings
can be performed in an arbitrary order. Nevertheless as we shall see this
construction plays an important role in composing nontrivial chains.

The previous example also demonstrates that it is worth searching the
conditions which
will guarantee that each member of a sequence of elements of
the type \eiphi \
is the solution of the equations (\ref{gentwist}) for coproducts defined by all
the previous twists of this sequence.

One of the obvious solutions to this problem can be formulated as follows:

\newtheorem{theor1}{Proposition}
\begin{theor1}\label{prechain}
Let ${\cal A}$ be a Hopf algebra, ${\cal B}$ and ${\cal C}$ be its
subalgebras such that they are carriers for twists ${\cal F_B}$ and  ${\cal
F_C}$ respectively. Let ${\cal F_B}$ commute with $\Delta {\cal C}$.
Then ${\cal C}$ is stable with respect to ${\cal F_B}$,
${\cal F_C}$ is a twisting element for ${\cal A}_{{\cal F}_{\cal B}}$ and
the composition
\be
{\cal F_C}{\cal F_B}
\ee
is a twisting element for ${\cal A}$.
\end{theor1}
In the previous example ${\cal B}$ and ${\cal C}$ were the Heisenberg
subalgebras in $sl(N)$ intersecting by the element $E_{1N}$. The other trivial
case is when ${\cal A}$ contains the direct sum ${\cal B} \oplus {\cal C}$.
In the next section we shall study some nontrivial cases typical for the
universal enveloping classical Lie algebras.

\sect{Chains}

For the classical Lie algebras there exists the possibility to construct
sequences of carrier subalgebras systematically.
\begin{theor1}
Let ${\cal A}$ be a Hopf algebra and
\be
\label{gentower}
{\cal A}_p \subset {\cal A}_{p-1} \subset \ldots \subset {\cal A}_{1} \subset
{\cal A}_{0} \equiv {\cal A}
\ee
a sequence of Hopf subalgebras such that
\be
\label{inject}
 {\cal B}_{k} \subset {\cal A}_{k} , \quad k=0,1, \ldots , p,
\ee
are the carrier subalgebras for twisting elements ${\cal F}_{{\cal B}_k}$.
Let  ${\cal F}_{{\cal B}_k}$ commute with
$ \Delta {\cal A}_{k+1}$:
\be
\fabk \Delta {\cal A}_{k+1} = \Delta {\cal A}_{k+1}  \fabk
\ee
 Then for any $k=0,1, \ldots, p$ the composition
\be
\chfabk \equiv
{\cal F}_{{\cal B}_k}{\cal F}_{{\cal B}_{k-1}} \ldots {\cal F}_{{\cal B}_0}
\ee
is a twisting element for ${\cal A}$.
\end{theor1}

Now we shall show how this scheme can be realized for the universal enveloping
algebras $U(\g)$ for classical Lie algebras $ \g $ ( $U(\g)$ is considered
here as a Hopf algebra with primitive comultiplication of generators).
The construction will
be similar for the classical series $A_n$, $B_n$ and $D_n$. In
the case of simplectic algebras $C_n$ the chain would not be completely
proper and we shall treat this situation separately.

Let us consider the following sequences of Hopf algebras:
\be
\label{sltower}
 \begin{array}{l}
 U(sl(N)) \supset U(sl(N-2)) \supset \ldots \supset U(sl(N-2k)) \supset
 \ldots \\
 \end{array}
 \,\, {\rm for } \,\, A_{N-1}
\ee
\be
\label{soetower}
  \begin{array}{l}
U(so(2N)) \supset U(so(2(N-2)) \supset \ldots \supset U(so(2(N-2k)) \supset
 \ldots \\
\end{array}
 \,\, {\rm for } \,\, D_{N}
\ee
\be
\label{sootower}
  \begin{array}{l}
U(so(2N+1)) \supset U(so(2(N-2)+1) \supset \ldots \supset U(so(2(N-2k)+1)
\supset \ldots  \\
\end{array}
 \,\, {\rm for } \,\, B_{N}
\ee
We want to show that for these sequences there exist the sets of maximal
${\cal ET}$'s with
the properties listed
in the Proposition 2. In each element \cak of the sequences
let us choose the ``initial"
root \lamnk . All the roots are equivalent in
$A$ and $D$ series, but in the series $B$ one of the long roots must be chosen
(this will be justified later). For definiteness we fix the following choice
(all the roots are written in the
$e$-basis):
\begin{equation}
\label{inroots}
 \lamnk =
\left\{
\begin{array}{lcl}
 e_1 - e_2  & {\rm for}
& sl(N-2k),\\
 e_1 + e_2
& {\rm for} & so(2(N-2k)),\\
 e_1 + e_2
& {\rm for} & so(2(N-2k)+1),
\end{array}
\right.
\end{equation}

Consider the set of
roots orthogonal to \lamnk. They form the subsystems for the following
subalgebras in \cak:
\be
\label{inj1}
sl(M-2) \subset sl(M),
\ee
\be
\label{inj2}
so(2M-4) \oplus sl(2) \subset so(2M),
\ee
\be
\label{inj3}
so(2M-3) \oplus sl(2) \subset so(2M+1).
\ee
Notice that in all the cases (\ref{inj1})-(\ref{inj3})
the corresponding universal enveloping algebras contain \cakp.

For each \cak let us form the set \kpi of roots $\lambda$ that are the
constituent for the initial root \lamnk, i.e.,
\be
\label{kpi}
\kpi = \left\{ \lambda',\lambda''|\lambda'+\lambda'' = \lamnk ; \quad
\lambda' + \lamnk,\lambda' + \lamnk \not\in \Lambda_{\cal A} \right\} .
\ee
For each element $\lambda ' \in \kpi$ one can indicate
the root $\lambda '' \in \kpi$ such that $\lambda ' + \lambda '' = \lamnk$.
Let us consider the ordered pairs $(\lambda',\lambda'')$ and
decompose \kpi according to its main property
\be
\begin{array}{l}
\kpi =   \pi_k' \, \cup \,  \pi_k'' ,\\
 \pi_k' = \{ \lambda' \}, \quad  \pi_k''= \{ \lambda'' \}.
\end{array}
\ee
For the sequences we are
dealing with these sets are:
\begin{equation}
\label{kpiroots}
\{ \lambda'_l,\lambda''_l \} =
\left\{
\begin{array}{lcl}
\left\{ \left\{(e_1 - e_l) \right\},
\left\{(e_l -e_2) \right\} \right\}_{l=3,4, \ldots ,M} & {\rm for}
& sl(M),\\
\left\{ \left\{(e_1 \pm e_l) \right\}, \left\{(e_2 \pm e_l) \right\}
\right\}_{l=3,4,
\ldots ,M}
& {\rm for}
& so(2M),\\
\left\{ \left\{e_1, (e_1 \pm e_l) \right\}, \left\{ e_2, (e_2 \pm e_l)
\right\} \right\}_{l=3,4,
\ldots ,M}
& {\rm for} & so(2M+1),
\end{array}
\right.
\end{equation}

The important observation is that the generators $L_{\lambda'}$ and
$L_{\lambda''}$ for $\lambda' \in  \pi_k'$
 and $\lambda'' \in  \pi_k''$ form the bases for the spaces of conjugate
defining representations of the subalgebras $\cakp
\subset \cak$  (with respect to the adjoint action). These subrepresentations
are
\be
\label{slreps}
\left\{
\left( \underline{M-2} \right), \left( \underline{M-2} \right)^*
\right\} \quad {\rm for} \quad U(sl(M-2)) \subset U(sl(M)),
\ee
\be
\label{soereps}
\left\{
\left( \underline{2(M-2)} \right), \left( \underline{2(M-2)} \right)
\right\} \quad {\rm for} \quad U(so(2(M-2))) \subset U(so(2M)),
\ee
\be
\label{sooreps}
\left\{
\left( \underline{2M-3)} \right), \left( \underline{2M-3)} \right)
\right\} \quad {\rm for} \quad U(so(2M-3)) \subset U(so(2M+1)).
\ee
Notice that any generator $L_{\lambda}$ ( $\lambda \in \kpi$) commutes with
$L_{\scriptsize{\lamnk}}$ and with all
the other elements $\left\{ L_{\mu} | \mu \in \kpi \right\}$ except its
counterpart -- the generator $L_{\scriptsize{\lamnk - \lambda}}$.
Together with $L_{\lambda}$
we shall consider the Cartan generator \hlamnk dual
to the initial root (its
functional $(\hlamnk)^*$ is proportional to \lamnk) . To simplify
the  expressions we shall use the fact that
in any classical Lie algebra there
exists a basis where the structure constants for the generators
$\left\{ L_{\lambda}, H_{\scriptsize{\lamnk}}|\lambda
\in \kpi  \right\} $ can be
normalized to form the following compositions:
\begin{equation}
\label{algel}
\begin{array}{lcl}
 \left[ \hlamnk, L_{\lambda'} \right] & = & \frac{1}{2} L_{\lambda'},\\[2mm]
 \left[ \hlamnk, L_{\scriptsize{\lamnk} - \lambda' } \right]
 & = & \frac{1}{2}
                       L_{\scriptsize{\lamnk} - \lambda' },\\[2mm]
 \left[ L_{\lambda'}, L_{\scriptsize{\lamnk} - \lambda'} \right] & =
 & L_{\scriptsize{\lamnk}},
\end{array}
\begin{array}{lcl}
 \left[  L_{\scriptsize{\lamnk}}, L_{\lambda'} \right] & = &
\left[  L_{\scriptsize{\lamnk}}, L_{\scriptsize{\lamnk} - \lambda'} \right]
 = 0,
 \\[2mm]
 \left[ \hlamnk, L_{\scriptsize{\lamnk}} \right]
 & = &  L_{\scriptsize{\lamnk}},\\[2mm]
\lambda' \in \pi_k', & & \lamnk - \lambda' \in \pi_k''.
\end{array}
\end{equation}
In the example considered in Section 7 we present the explicit
realizations for the generators of
\cak  that fit the relations above.

The relations (\ref{algel}) show that for every
triple of roots $\left\{ \lambda' , \lamnk - \lambda', \lamnk \right\}$ we
have the subalgebra $\L_{\lambda'}(\alpha, \beta)$ with $\alpha = \beta =
\frac{1}{2}$  (see (\ref{l-norm})).
The set of generators
\be
\label{bk}
 \left\{ L_{\lambda |_{ \lambda \in \kpi}}, L_{\scriptsize{\lamnk}}, \hlamnk
 \right\}
\ee
define a subalgebra $\cbk \subset \cak$.

Let us perform in \cak the Jordanian twist
\be
\jkphi = \exp \{\hlamnk \otimes \signk \}
 \ee
with $\signk = \ln(1 + L_{\scriptsize{\lamnk}} )$.
In the twisted algebra $(\cak)_{\scriptsize{\jkphi}}$ the subalgebras
$\left\{ \L^{\lambda'}(1/2,1/2)\right.$
 $| \left. \lambda' \in \pi'_k \right\}$ described
above obviously obey the conditions of the
Proposition \ref{prechain} and the corresponding sequence of twists
\be
  \ekphi=\prod_{\lambda' \in \pi'_k}\ephi_{\lambda'} = \prod_{\lambda' \in
\pi'_k}
\exp \{  L_{\lambda'}
\otimes
 L_{\scriptsize{\lamnk} - \lambda'} e^{-\frac{1}{2} \signk } \}
\ee
 can be performed in it. This gives for each \cak the following ${\cal ET}$
 element:
\be
\label{fabk}
\fabk =  \ekphi \jkphi = \left( \prod_{\lambda' \in \pi'_k}\ephi_{\lambda'}
\right) \jkphi.
\ee

The sets of algebras \cak presented in (\ref{sltower}),(\ref{soetower})
and (\ref{sootower}) together with their subalgebras \cbk
(defined by the basic
families (\ref{bk})) form the correlated sequences of
subalgebras that satisfy
the conditions of the Proposition 2.
To prove this
let us consider the adjoint representation ${\rm ad}(\caz) \equiv d_0$
and its restrictions to the subalgebras \cak
: $d_k = {\rm ad}(\caz)_{|\scriptsize{\cak}}$. The
space of \cbk is invariant with respect to $d_{k+1}$. It contains the
subspaces of two trivial subrepresentations (generated by \lamnk and by
\hlamnk). This means that the ${\cal JT}$ factor \jkphi commutes
with the algebra $\Delta (\cakp) \subset \cakp \otimes \cakp $.
The other two invariant subspaces in \cbk refer
to the fundamental representations of
\cakp indicated in (\ref{slreps}), (\ref{soereps}) and (\ref{sooreps}).
Due to the commutation relations in \cbk the element $\ln \ekphi$ can be
written as
\be
\left( \sum_{\lambda' \in \pi'_k  } L_{\lambda'} \otimes
 L_{\scriptsize{\lamnk} - \lambda'} \right)  e^{-\frac{1}{2} \signk}
\ee
With the ordered pairs of roots as in (\ref{kpiroots}) this expression is
$d_{k+1}$-invariant (the converted conjugate bases for representations
modulo the scalar factor).
We have arrived at the conclusion that the sets of subalgebras
\cak ((\ref{sltower}), (\ref{soetower}) and (\ref{sootower})) and \cbk
(defined by  (\ref{kpi}), (\ref{kpiroots}) and (\ref{bk})) with the twisting
elements \fabk (\ref{fabk}) satisfy the conditions of the Proposition 2.
Thus for any classical simple Lie algebra of the series $A$, $B$ and $D$
the chains of twists
$ \chfabk \equiv {\cal F}_{{\cal B}_k}{\cal F}_{{\cal B}_{k-1}} \ldots
{\cal F}_{{\cal B}_0}  \, \, (k=0,1, \ldots, p) $ exist.

The twisting element for a chain can be written explicitely as
\be
\label{chain}
\begin{array}{l}
\chfabk =
\prod_{\lambda' \in
\pi'_k  }
\left( \exp \{  L_{\lambda'}
\otimes
 L_{\scriptsize{\lamnk} - \lambda'} e^{-\frac{1}{2} \signk } \} \right)
              \cdot \exp \{\hlamnk \otimes \signk \} \,\cdot \\
\prod_{\lambda' \in
\kmpip}
\left( \exp \{  L_{\lambda'}
\otimes
 L_{\scriptsize{\lamnkm} - \lambda'} e^{-\frac{1}{2} \signkm } \} \right)
               \cdot \exp \{\hlamnkm \otimes \signkm \} \, \cdot \\
\ldots \\
\prod_{\lambda' \in
\zpip}
\left( \exp \{  L_{\lambda'}
\otimes
 L_{\scriptsize{\lamnz} - \lambda'} e^{-\frac{1}{2} \signz } \} \right)
                \cdot \exp \{\hlamnz \otimes \signz \}

\end{array}
\ee
Any number of exponential factors can be cut out from the left. The remaining
part always conserves the property of the twisting element for the
corresponding classical Lie algebra. When on the left-hand side one has a
product of extensions that is not complete (not all  $\lambda' \in
\pi'_k  $ are used) ,
\be
\faethe \chfabkm =
\prod_{\lambda' \in \,
\Theta \,\subset \pi'_k  }
\left( \exp \{  L_{\lambda'}
\otimes
 L_{\scriptsize{\lamnk} - \lambda'} e^{-\frac{1}{2} \signk } \} \right)
              \cdot \exp \{\hlamnk \otimes \signk \} \cdot
{\cal F}_{{\cal B}_{k-1}}\cdot \ldots \cdot {\cal F}_{{\cal B}_0},
\ee
the subalgebra \cakp will be twisted nontrivially by such an element. In
this case the twisting deformation with \fabkp (of the (\ref{fabk}) type)
cannot be applied to ${\cal A}_{\scriptsize{\faethe \chfabk}}$. The necessary
{\sl primitivization of generators in
\cakp is regained when the product of extensions is complete} and
forms an invariant of the representation
$d_{k+1}$. We call this the ``matreshka" effect.

Quantizations  ${\cal A}_{\scriptsize{\chfabp}}$ of
classical Lie algebras produce the
chains of ${\cal R}$-matrices:
\be
\label{chrmat}
\cherbp = ({\cal F}_{{\cal B}_{p}})_{21}
( {\cal F}_{{\cal B}_{p-1}})_{21} \ldots
( {\cal F}_{{\cal B}_{0}})_{21}
{\cal F}_{{\cal B}_{0}}^{-1} \ldots
{\cal F}_{{\cal B}_{p-1}}^{-1}
{\cal F}_{{\cal B}_{p}}^{-1}.
\ee
The explicit expressions in terms of generators can be obtained substituting the
elements $L_{\lambda}$ and $\hlamnk$ in (\ref{chain}) by the corresponding
generators according to the prescription of roots in (\ref{inroots}) and
(\ref{kpiroots}).

If the deformation parameter is introduced (as in (\ref{twist-sl(N)})) the chains
of classical $r$-matrices can be extracted from (\ref{chrmat}):
\be
\label{chcrmat}
\chercbp = \sum_{k=0,1,\ldots,p} \left( \hlamnk \wedge L_{\lambda^k_0} +
\sum_{\lambda' \in \kpi} L_{\lambda'} \wedge
L_{\lambda^k_0 -\lambda'} \right).
\ee

With the obvious modifications of factors  (summands) the sequences of ${\cal
R}$-matrices (classical $r$-matrices) for  incomplete chains of twists can
also be written.

\sect{Improper chains. Simplectic algebras }

Imposing additional conditions on the internal structure of the Hopf algebras
involved one can minimize the algebra \caz on which the chain is based to the
universal enveloping algebra $\cazcar \equiv U(\gzcar) $ of the carrier
of the chain. This happens, for example, when \gz in $\caz = U(\gz)$ is a
sequence of  semidirect sums and every \cbk is a \cbkp-module with respect to the
adjoint action (in \gz),
\be
\begin{array}{c}
\label{addcond}
\,\gzcar   = \gz = {\cal B}_p \vdash ({\cal B}_{p-1} \vdash (
\cdots \vdash {\cal B}_0)\cdots),\\[2mm]
\, \left[ \cbkp,\cbk \right] \subset \cbk.
\end{array}
\ee
In this case one can define the
subalgebras $\cakcar   $ as
\be
\label{mincak}
U(\gkcar)   =U \left( {\cal B}_p \vdash ({\cal B}_{p-1} \vdash (
\cdots \vdash {\cal B}_k)\cdots) \right).
\ee
In the sequences of classical algebras that we considered in (\ref{sltower}),
(\ref{soetower}) and (\ref{sootower}) the conditions (\ref{addcond}) are
fulfilled (with \cbk defined by (\ref{bk})). One can rewrite the sequences
(\ref{gentower}) for  the classical Lie algebras so that the
elements \cak will be substituted by $\cakcar \equiv U(\gkcar) $
defined by (\ref{mincak}) and $\gzcar  $ will be the carrier of
\chfabp. There rests some freedom in choosing the initial root.
Using it one can, in particular, place the carrier of the chain
in the Borel subalgebra of the corresponding classical Lie algebra.
For example, in the case of $sl(N)$ the carrier subalgebra of the
full chain of the type (\ref{chain}) can be arranged to
contain all the generators with the positive root vectors and
a part of the Cartan subalgebra (spanned by $H_{1,N}, H_{2,N-1}, \ldots$).

For simple Lie algebras the chain carrier subalgebra covers only a
proper subspace of an algebra. The chains \chfabp that we described
in the previous section are maximal in the sence that $\cazcar  $ is
a maximal Frobenius subalgebra in the corresponding classical Lie algebra
that can be composed from the subalgebras of the type (\ref{algel})
(that is, using \ekphi \ and \jkphi \ as elementary blocks). These
chains are also specific for the simple algebras we are dealing with.
In each of the three cases ((\ref{sltower}), (\ref{soetower}) and
(\ref{sootower})) the individual properties of the root system are
used to form a chain.

The universal enveloping algebras for other simple Lie algebras
(the series $C_N$ and the exeptional algebras) do not refer to
the class of algebras conserving symmetric forms
(over a field) and cannot be supplied by a specific chain of extended
twists. Nevertheless, the quantization by a chain
can be performed in these algebras using
the classical subalgebras of the series $A,B$ and $D$ contained in them.
For example, due to the  inclusion $sl(N) \subset sp(N)$ the chain
specific to $U(sl(N))$ can be used to quantize $U(sp(N))$. Such chains
can be called improper.

Now we shall study the universal enveloping algebras
for simplectic simple Lie algebras
(the series $C_N$) where the maximal chain appears to be improper. (It
exploits almost exclusively the $A_{N-1}$ subalgebra in $C_N$.) In the
$e$-basis the $sp(N)$ roots can be fixed as follows
\be
\label{sproots}
\Lambda_{sp(N)} =
\left\{
\begin{array}{lcl}
e_i - e_j & {\rm for} & i \neq j \\
\pm(e_i + e_j)  & {\rm for} & i \leq j
\end{array}
\right\}
\quad
i,j = 1,2, \ldots , N.
\ee
Whatever root will be chosen as the initial one the extensions will contain
generators whose roots will have the nonzero projections on the
$sp(N-2)$ root system.

Note that if we fix a short root to be the initial,
for example $\lambda_0^0 = e_i -
e_j$, there will be  pairs of constituent roots that do not satisfy the
conditions (\ref{kpi}). The generators corresponding to $\lambda' = -2e_j$,
$\lambda'' = e_i + e_j$ and
$\lambda^0_0$ do not form a subalgebra of $\L(\alpha,\beta)$-type.
Thus we are to consider
the subalgebra ${\cal A}_1 = U(sp(N-2))$. The generators corresponding to
$\pi_k'$ and $\pi_k''$ (\ref{kpi}),
$$
\left\{ e_i \pm e_l \right\} \quad {\rm and} \quad \left\{ -e_j \pm e_l
\right\} \quad \quad l=3,\dots ,N,
$$
form the bases for the defining representations of $sp(N-2)$. The simplectic
invariant
\be
\label{simp}
\sum_l \left(  L_{e_i + e_l} \otimes  L_{-e_j - e_l} - L_{e_i - e_l} \otimes
L_{-e_j + e_l}\right)
\ee
does not correlate with the ${\cal ET}$ (\ref{fractions}).
Otherwise one can check that the extensions
based on linear combinations of the
type (\ref{simp}) (with the coefficients in
${\bf C}\left[\left[\sigma \right]\right]$) violate the twist equation
(\ref{gentwist}).

We can diminish the subalgebra ${\cal A}_1$ and put ${\cal A}_1
= U(sl(N-2))$. In this case the summands $\sum_l  L_{e_i + e_l}
\otimes  L_{-e_j - e_l}$ and $\sum_l  L_{e_i - e_l} \otimes L_{-e_j + e_l}$
in (\ref{simp}) will be separately invariant with respect to ${\cal A}_1$
and both will match with the sequences of extensions  (\ref{fractions}).
In such a way we can proceed constructing the chain of extended twists
for $U(sp(N))$ but this will be specific for $A_n$ rather than
for $C_n$ root system (except that the long root can be chosen
to be the first initial root).

\sect{Example. Maximal chain for $U(so(9))$}
To illustrate the properties of chains we apply the algorithm presented
in Sections 5 and 6 to construct a chain of ${\cal ET}$'s for the
algebra $U(so(9))$.

In this case the sequence (\ref{sootower}) consists of two elements:
\be
\label{so9seq}
\cao \supset \caz = so(9) \supset so(5)
\ee
with the initial roots
$$
\lambda^0_0 = e_1 + e_2 , \quad \quad \lambda^1_0 = e_3 + e_4
$$
and the corresponding sets of constituent roots
$$
\begin{array}{lcl}
\pi'_0 & = \{ \lambda^{0 \prime} \} =
& \{ e_1, e_1 \pm e_3, e_1 \pm e_4  \}\\[2mm]
\pi''_0 & = \{ \lambda^{0 \prime \prime} \} =
& \{ e_2, e_2 \pm e_3, e_2 \pm e_4  \}\\[2mm]
\pi'_1 & = \{ \lambda^{1 \prime} \} = & \{  e_3  \}\\[2mm]
\pi''_1 & = \{ \lambda^{1 \prime \prime} \} = & \{  e_4  \}\\
\end{array}
$$
The roots $\pi'_0$ and $\pi''_0$ form the weight diagrams for the vector
representations of $\cao= so(5)$.

Together with the Cartan generators $H_{\scriptsize{\lamnz}}$
, $H_{\scriptsize{\lamno}}$ the
basic elements $\{ E_{\lambda} |\lambda \in \zpi \cup \opi \}$ and
$E_{\scriptsize{\lamnz}}, E_{\scriptsize{\lamno}}$ form the 16-dimensional
subalgebra $\gzcar \subset \gz = so(9) $. It has the
structure of a semidirect sum $\gzcar \approx \cbo \vdash \cbz$.
This means that studying this chain we can restrict ourselves to
the subalgebra $U(\gzcar)$.

The maximal chain for the sequence (\ref{so9seq}) has the following
structure
\be
\label{so9ch}
\chfabo = \eophi \jophi \ezphi \jzphi = \Phi_{{\cal E}\lambda_3}
\jophi ( \prod_{\lambda' \in \pi'_0} \Phi_{{\cal E}\lambda'} )
\jzphi .
\ee

The generators of $\gz$ can be expressed in terms of the antisymmetric
Okubo matrices  $M_{ik}$ :
\be
\begin{array}{l}
{\bf L}^{12} = \left\{
\begin{array}{l}
H_{12} = (-i/2)(M_{12} + M_{34}), \\
E_{1} = M_{29} -i M_{19}, \quad \quad \quad \quad \quad
E_{2} = M_{49} -i M_{39},\\
E_{1+2} = -M_{24} + i M_{23} +iM_{14} + M_{13} ,  \\
\end{array}
\right.\\[4mm]
\begin{array}{l}
E_{1+3}  = -M_{26} + i M_{25} +iM_{16} + M_{15} , \\
E_{1+4}  = -M_{28} + i M_{27} +iM_{18} + M_{17} , \\
E_{2+3}  = -M_{46} + i M_{45} +iM_{36} + M_{35} , \\
E_{2+4}  = -M_{48} + i M_{47} +iM_{38} + M_{37} , \\
E_{1-3} = -M_{26} - i M_{25} +iM_{16} - M_{15} , \\
E_{1-4} = -M_{28} - i M_{27} +iM_{18} - M_{17} , \\
E_{2-3} = -M_{46} - i M_{45} +iM_{36} - M_{35} , \\
E_{2-4} = -M_{48} - i M_{47} +iM_{38} - M_{37} , \\
\end{array}\\
\L^{34}  = \left\{
\begin{array}{l}
H_{34} = (-i/2)(M_{56} + M_{78}),  \\
E_{3} = M_{69} -i M_{59}, \quad \quad \quad \quad \quad
E_{4} = M_{89} -i M_{79},\\
E_{3+4}  = -M_{68} + i M_{67} +iM_{58} + M_{57} , \\
\end{array}
\right.\\[4mm]
\end{array}
\ee
Here the lower indices of raising generators $E$ indicate the
corresponding $so(9)$-roots.
The set of generators $\L^{12} $  ( $\L^{34} $) forms the
4-dimensional subalgebra of the type ${\bf L}(1/2, 1/2)$ with
$E=E_{1+2} $ ($E=E_{3+4} $).

The explicit expressions for the main factors of the chain in this basis
are as follows:
\be
\label{sofactors}
\begin{array}{l}
\jzphi = \exp(H_{12} \otimes \sigma_{12}), \quad
\jophi = \exp(H_{34} \otimes \sigma_{34}),\\
\ezphi = \exp (E_1 \otimes E_2 +
1/2 (E_{1-3} \otimes E_{2+3}  + E_{1+3}  \otimes E_{2-3}\\
\rule{3cm}{0cm} + E_{1-4} \otimes E_{2+4}  + E_{1+4}  \otimes E_{2-4} )
(1 \otimes e^{-\frac{1}{2}\sigma_{12}})),\\
\eophi = \exp(E_3 \otimes E_4 e^{-\frac{1}{2}\sigma_{34}})
\end{array}
\ee
with
$$
\begin{array}{l}
\sigma_{12}=\sigma_0^0=\ln(1+E_{1+2} ),\\
\sigma_{34}=\sigma_0^1=\ln(1+E_{3+4} ).\\
\end{array}
$$
After the first Jordanian twisting,
$$
  \cazcar \stackrel{\jzphi}{\longrightarrow} (\cazcar)_{J_0} ,
$$
the subalgebra
$$ \L^{34}  = {\cal B}_1 \subset {\cal A} $$
remains primitive. The carrier subalgebra for \jzphi \ acquires
the coproducts
\be
\label{bor12}
\begin{array}{l}
\Delta_{J_0} (H_{12}) = H_{12} \otimes e^{-\sigma_{12}}
                      + 1 \otimes H_{12},\\
\Delta_{J_0} (E_{1+2} ) = E_{1+2}  \otimes e^{\sigma_{12}}
                        + 1 \otimes E_{1+2},\\
\end{array}
\ee
The coproducts for the remaining generators of ${\cal B}_0$ are of the form
\be
\Delta_{J_0} (E) = E \otimes e^{\frac{1}{2}\sigma_{12}} + 1 \otimes E.
\ee

Among the exponential factors $\Phi_{{\cal E}\lambda'}$ of the
extension \ezphi \ (see(\ref{so9ch})) there is
one ($\Phi_{{\cal E}\lambda_1}$)
that does not touch the subalgebra ${\bf L}^{34}$. Each of the rest
$\{ \Phi_{{\cal E}\lambda'|\lambda' = e_1 \pm e_3 ,
e_1 \pm e_4} \}$ being applied separately produces a nontrivial
deformation of  ${\bf L}^{34}$. These extensions can be combined
to form the $so(5)$-invariant (see (\ref{sofactors})). In this case, i. e.
after the twisting
$$
 (\cazcar)_{J_0} \stackrel{\ezphi}{\longrightarrow}
 (\cazcar)_{{\cal E}_{0}J_0} ,
$$
the primitivity of generators in ${\bf L}^{34}$ is restored.
The coproducts for the generators of $(\cbz)_{{\cal E}_0 J_0}$ are
deformed according to the general rule (see Section 2),
\be
\begin{array}{lcl}
\Delta_{{\cal E}_0 J_0}(E_{\lambda'})
&=&
E_{\lambda'} \otimes e^{-\frac{1}{2}\sigma_{12}} + 1 \otimes E_{\lambda'},
\\
\Delta_{{\cal E}_0 J_0}(E_{\scriptsize{\lamnz} -\lambda'})
&=&
E_{\scriptsize{\lamnz} -\lambda'} \otimes
e^{\frac{1}{2}\sigma_{12}}
+ e^{\sigma_{12}} \otimes E_{\scriptsize{\lamnz} -\lambda'},
\\
\Delta_{{\cal E}_0 J_0}(E_{\scriptsize{\lamnz} })
&=&
E_{\scriptsize{\lamnz} } \otimes
e^{\sigma_{12}}
+ 1 \otimes E_{\scriptsize{\lamnz} },
\\
\Delta_{{\cal E}_0 J_0}(H_{12})
&=&
H_{12 } \otimes
e^{-\sigma_{12}} + 1 \otimes H_{12}
- E_1 \otimes E_2 e^{-\frac{3}{2}\sigma_{12}}\\
&&
-\frac{1}{2}\sum_{\lambda' = e_1 \pm e_3, e_1 \pm e_4}
E_{\lambda'} \otimes E_{\scriptsize{\lamnz} - \lambda'}
e^{-\frac{3}{2}\sigma_{12}}.
\\
\end{array}
\ee

As a result of the ``matreshka" effect the second Jordanian twist can be
applied to $(\cazcar)_{{\cal E}_0 J_0}$
$$
 (\cazcar)_{{\cal E}_0 J_0} \stackrel{\jophi}{\longrightarrow}
 (\cazcar)_{J_1{\cal E}_0 J_0} .
$$
This leads to the following deformations:
\begin{itemize}
\item the subalgebra \cbo acquires the well known twisted form with the
defining coproducts
\be
\begin{array}{l}
\Delta_{J_1{\cal E}_0 J_0} (H_{34}) = H_{34} \otimes
e^{-\sigma_{34}} + 1 \otimes H_{34} ,\\
\Delta_{J_1{\cal E}_0 J_0} (E_{3+4} ) = E_{3+4}  \otimes
e^{\sigma_{34}} + 1 \otimes E_{3+4} ,\\
\Delta_{J_1{\cal E}_0 J_0} (E_k) = E_k \otimes
e^{\frac{1}{2}\sigma_{34}} + 1 \otimes E_k, \quad \quad k=3,4;
\end{array}
\ee
\item the subalgebra $({\bf L}^{12})_{{\cal E}_{0} J_0}$ rests untouched
$$
({\bf L}^{12})_{{\cal E}_{0} J_0} = ({\bf L}^{12})_{J_1{\cal E}_{0}J_0};
$$
\item for each $\{ {\lambda = e_i \pm e_k} | i=1,2; k=3,4 \}$
the following substitution is performed in the coproducts for
the generators $E_{\lambda}$
$$
E_{\lambda} \otimes f(\sigma_{12})
\quad \longrightarrow
\quad E_{\lambda} \otimes e^{\pm\frac{1}{2}\sigma_{34}}f(\sigma_{12});
$$
\item in
$\Delta_{J_1 {\cal E}_0 J_0}(E_{\lambda})$ for each
$\{ \lambda= e_i - e_k \ | \ i=1,2; k=3,4 \}$
the additional term appears,
$$
(-1)^{k+1}H_{34}e^{(i-1)\sigma_{12}} \otimes
E_{e_i + e_{\overline{k}}} e^{-\sigma_{34}}
$$
(here $\overline{3}=4, \quad \overline{4}=3$);
\item for the Cartan generator $H_{12}$ the coproduct becomes
$$
\begin{array}{lcl}
\Delta_{J_1 {\cal E}_0 J_0}(H_{12})
& = &
H_{12} \otimes e^{-\sigma_{12}} + 1 \otimes H_{12}
- E_1 \otimes E_2 e^{- \frac{3}{2}\sigma_{12}}\\
&&
-\frac{1}{2} E_{1+3}  \otimes E_{2-3}
e^{\frac{1}{2}\sigma_{34}- \frac{3}{2}\sigma_{12}}
- H_{34} E_{1+3}  \otimes E_{2+4}
e^{-\frac{1}{2}\sigma_{34}- \frac{3}{2}\sigma_{12}}\\
&&
-\frac{1}{2} E_{1+4}  \otimes E_{2-4}
e^{\frac{1}{2}\sigma_{34}- \frac{3}{2}\sigma_{12}}
+H_{34} E_{1+4}  \otimes E_{2+3}
e^{-\frac{1}{2}\sigma_{34}- \frac{3}{2}\sigma_{12}}\\
&&
-\frac{1}{2} E_{1-4} \otimes E_{2+4}
e^{-\frac{1}{2}\sigma_{34}- \frac{3}{2}\sigma_{12}}
-\frac{1}{2} E_{1-3} \otimes E_{2+3}
e^{-\frac{1}{2}\sigma_{34}- \frac{3}{2}\sigma_{12}}.\\
\end{array}
$$
\end{itemize}
The last twisting (that completes the chain ${\cal B}_{1 \prec 0}$),
\be
\label{lasttw}
 (\cazcar)_{J_1{\cal E}_0 J_0} \stackrel{\eophi}{\longrightarrow}
 (\cazcar)_{{\cal B}_{1 \prec 0}},
\ee
does not change the coproducts for the
generators $\{ E_{i+k} \ | \ i=1,2; \ k=3,4  \}$.
It produces the ordinary transformation for $\L^{34} $,
$$
\begin{array}{lcl}
\Delta_{{\cal B}_{1 \prec 0}}(H_{34}) &= &
 \Delta_{J_1 {\cal E}_0 J_0}(H_{34})
+ E_3 \otimes E_4 e^{-\frac{3}{2}\sigma_{34}}\\[1mm]
\Delta_{{\cal B}_{1 \prec 0}}(E_{3})  &= &
 E_{3} \otimes e^{-\frac{1}{2}\sigma_{34}}
 + 1 \otimes E_{3}\\[1mm]
\Delta_{{\cal B}_{1 \prec 0}}(E_{4})  &= &
 E_{4} \otimes e^{\frac{1}{2}\sigma_{34}}
 + e^{\sigma_{34}} \otimes E_{4}.\\
\end{array}
$$
The generators
$E_1, E_2, E_{i-k}$ and $H_{12} $
are  nontrivially twisted by the transformation (\ref{lasttw}),
$$
\begin{array}{lcl}
\Delta_{{\cal B}_{1 \prec 0}}(E_1) & = &
 \Delta_{J_1{\cal E}_0 J_0}(E_1)
- E_{1+3}  \otimes E_4 e^{-\frac{1}{2}\sigma_{34}-\frac{1}{2}\sigma_{12}}
- E_{3} \otimes E_{1+4}  e^{-\frac{1}{2}\sigma_{34}},\\[1mm]
\Delta_{{\cal B}_{1 \prec 0}}(E_2) & = &
\Delta_{J_1{\cal E}_0 J_0}(E_2)
- E_{2+3}  \otimes E_{4} e^{-\frac{1}{2}\sigma_{34}+\frac{1}{2}\sigma_{12}}
- E_{3}e^{\sigma_{12}} \otimes E_{2+4}  e^{-\frac{1}{2}\sigma_{34}},\\[1mm]
\Delta_{{\cal B}_{1 \prec 0}}(E_{1-3}) & = &
 \Delta_{J_1{\cal E}_0 J_0}(E_{1-3})
+ 2E_{1} \otimes E_4 e^{-\sigma_{34}-\frac{1}{2}\sigma_{12}}\\[1mm]
& &
 - E_{1+3}  \otimes E_{4}^2
e^{-\frac{1}{2}\sigma_{12}-\frac{3}{2}\sigma_{34}}
-2E_{3} \otimes E_4 E_{1+4}  e^{-\frac{3}{2}\sigma_{34}}
,\\[1mm]
\Delta_{{\cal B}_{1 \prec 0}}(E_{2-3}) & = &
 \Delta_{J_1{\cal E}_0 J_0}(E_{2-3})
+ 2E_{2} \otimes E_4 e^{-\sigma_{34}+\frac{1}{2}\sigma_{12}}\\[1mm]
& &
 - E_{2+3}  \otimes E_{4}^2
e^{\frac{1}{2}\sigma_{12}-\frac{3}{2}\sigma_{34}}
-2E_{3}e^{\sigma_{12}} \otimes E_4 E_{2+4}  e^{-\frac{3}{2}\sigma_{34}}
,\\[1mm]
\Delta_{{\cal B}_{1 \prec 0}}(E_{1-4}) & = &
 \Delta_{J_1{\cal E}_0 J_0}(E_{1-4})
+ 2E_{3} \otimes E_1 e^{-\frac{1}{2}\sigma_{34}}\\[1mm]
& &
 - E_{3}^2 \otimes E_{1+4}  e^{-\sigma_{34}}
+2E_{3} \otimes E_4 E_{1+3}  e^{-\frac{3}{2}\sigma_{34}}
,\\[1mm]
\Delta_{{\cal B}_{1 \prec 0}}(E_{2-4}) & = &
 \Delta_{J_1{\cal E}_0 J_0}(E_{2-4})
+ 2E_{3}e^{\sigma_{12}} \otimes E_2 e^{-\frac{1}{2}\sigma_{34}}\\[1mm]
& &
 - E_{3}^2 e^{\sigma_{12}}\otimes E_{2+4}  e^{-\sigma_{34}}
+2E_{3} e^{\sigma_{12}} \otimes E_4 E_{2+3}
e^{-\frac{3}{2}\sigma_{34}},\\[1mm]
\end{array}
$$
$$
\begin{array}{lcl}
\Delta_{{\cal B}_{1 \prec 0}}(H_{12}) & = &
 \Delta_{J_1{\cal E}_0 J_0}(H_{12})
+ (E_{1+3}  \otimes E_{2}E_{4} + E_{3}E_{1} \otimes E_{2+4})
 1 \otimes e^{-\frac{1}{2} \sigma_{34}-\frac{3}{2}\sigma_{12}}\\[1mm]
& &
 - (E_{1} + E_{1+4} E_{3}) \otimes E_{2+3}  E_4
 e^{-\sigma_{34}-\frac{3}{2}\sigma_{12}}
  + \frac{1}{2} E_{1+3}  \otimes E_{2+3}  (E_4)^2
 e^{-\frac{3}{2}(\sigma_{34}+\sigma_{12})}\\[1mm]
& &
 - E_{1+4} E_{3} \otimes E_{2}
 e^{-\frac{3}{2}\sigma_{12}}
  + \frac{1}{2} E_{1+4}  (E_{3})^2 \otimes E_{2+4}
 e^{-\frac{1}{2}\sigma_{34}-\frac{3}{2}\sigma_{12}}.
 \end{array}
$$

These relations complete the description of the twisted Hopf algebra
$(\cazcar)_{{\cal B}_{1 \prec 0}}$. Using the explicit expressions
(\ref{sofactors}) for the chain factors  one can reconstruct the
Hopf algebra $U(so(9))_{{\cal B}_{1 \prec 0}}$
containing $(\cazcar)_{{\cal B}_{1 \prec 0}}$. Both of them are
triangular with the universal element
$$
\cherbk = \left( \eophi \jophi \ezphi \jzphi \right)_{21}
\left( \eophi \jophi \ezphi \jzphi \right)^{-1}.
$$
The deformation parameter can be introduced so that the
classical $r$-matrix
$$
\begin{array}{lcl}
r_{{\cal B}_{1 \prec 0}} & = & H_{12} \wedge E_{1+2}  +
H_{34} \wedge E_{3+4}  + E_1 \wedge E_2 + E_3 \wedge E_4\\
& &
+ 1/2 (E_{1-3} \wedge E_{2+3}  + E_{1+3}  \wedge E_{2-3}
+ E_{1-4} \wedge E_{2+4}  + E_{1+4}  \wedge E_{2-4} ).\\
\end{array}
$$
determines the Lie-Poisson structure that was quantized explicitly
by the chain of twists \chfabk.
\sect{Conclusions}
Chains of twists provide a rich variety of new quantizations for
a certain class of Lie algebras described in Proposition 2.
As it was demonstrated in \cite{V-M-2} extended twists can be
accompanied by the special Reshetikhin twists which ``rotate" the roots
of the carrier subalgebras for the Jordanian factors. It is easy
to check that such ``rotations" can be applied also in the case of
chains. The corresponding additional factors
${\cal F}_{\cal R}=\exp \left( \left(\hlamnkm +
\theta (\hlamnkm)^{\perp} \right) \otimes \sigma^{k-1}_0  \right)$
(here $(\hlamnkm)^{\perp}$ is orthogonal to \hlamnk and \hlamnkm) can be
included in each \fabkm. It can be shown that though
the factor \ekmphi must be changed its invariance properties
with respect to \cbk can be conserved. In this context the chains are
flexible and their multiparametric versions can be easily constructed.

The deformation parameters can be introduced in chains by rescaling
the generators of the subalgebra \cbk. It must be stressed that each
\cbk can be rescaled separately with an independent variable $\xi_k$.
When all these rescaling factors are proportional to the deformation
parameter $\xi$, i.e. $\xi_k = \xi \eta_k $,
then in the classical limit the
parameters $\eta_k$ appear as the multipliers in the classical $r$-matrix
(compare with (\ref{chcrmat})):
$$
\chercbp = \sum_{k=0,1,\ldots,p} \eta_k \left( \hlamnk \wedge
L_{\lambda^k_0} +
\sum_{\lambda' \in \kpi} L_{\lambda'} \wedge
L_{\lambda^k_0 -\lambda'} \right).
$$

The mechanisms described above can be combined together both leading to the
multiparametric versions of chains.

One of the consequences of the Proposition 2 is that for a
large set of universal enveloping algebras (including $A, \, B$ and $D$
series of classical algebras) the
classical $r$-matrices of the type (\ref{chcrmat}) exist.
For the special case of $\g = sl(N)$ they were first presented in
\cite{GER}. As we have shown above they originate from the specific
properties of extended Jordanian twists -- the possibility to form
chains for the certain types of universal enveloping algebras.

\section*{Acknowledgments}

One of the authors (V. L.) would like to thank the DGICYT
of the Ministerio de Educaci\'on y Cultura de  Espa\~na for supporting
his sabbatical stay (grant SAB1995-0610). This  work has been partially
supported by DGES of the Ministerio de  Educaci\'on y Cultura of
Espa\~na under Project PB95-0719, the Junta de Castilla y Le\'on
(Espa\~na) and the  Russian Foundation for Fundamental Research
under the grants 97-01-01152 and 98-01-00310.


\end{document}